\documentclass[12pt,a4paper]{article}
\usepackage[active]{srcltx}
\usepackage{amsfonts}
\usepackage{amsmath}
\usepackage{graphicx,epsfig}
\usepackage{amssymb}
\usepackage{theorem,subfigure}
\usepackage[usenames,dvipsnames]{color} 
\usepackage{colortbl}
\definecolor{gray73}{RGB}{186,186,186}
\usepackage{graphicx}
\usepackage{array}
\usepackage{tikz}
\usetikzlibrary{matrix,arrows}
\usepackage{verbatim}
\usepackage{geometry}
\usepackage{multirow}

\usepackage{capt-of}
\usepackage{caption}
\DeclareCaptionLabelFormat{mylabelformat}{#1 #2:}
\newtheorem{theorem}{Theorem}[section]
\newtheorem{lemma}{Lemma}[section]
\newtheorem{remark}[theorem]{Remark}

\newtheorem{algorithm}[theorem]{Algorithm}

\numberwithin{equation}{section}

\newcommand{\C}{\mbox{$\mathbb C$}}
\newcommand{\R}{\mbox{$\mathbb R$}}

\newcommand{\supp}{\mbox{supp}\,}
\renewcommand{\mod}{\mbox{mod}\,}
\newcommand{\smod}{\mbox{\rm\scriptsize mod}\,}

\DeclareMathOperator{\e}{\mathrm e}
 at12pt

\title{A Sparse Fast Fourier Algorithm for Real Non-negative Vectors}
\author{Gerlind Plonka\footnote{University of G\"ottingen, Institute for Numerical and Applied Mathematics, Lotzestr.\ 16-18,  37083 G\"ottingen, Germany. Email: plonka@math.uni-goettingen.de} \qquad
Katrin Wannenwetsch\footnote{University of G\"ottingen, Institute for Numerical and Applied Mathematics, Lotzestr.\ 16-18,  37083 G\"ottingen, Germany. Email: k.wannenwetsch@math.uni-goettingen.de}}

\begin{document}
\maketitle

\abstract{In this paper we propose a new fast Fourier transform to recover a real non-negative signal ${\mathbf x} \in {\R}^{N}_{+}$ from its discrete Fourier transform $\widehat{\mathbf x} = {\bf F}_{N} {\mathbf x} \in {\C}^{N}$.
If the signal ${\mathbf x}$ appears to have a short support, i.e., vanishes outside a support interval of length 
$m < N$, then the algorithm has an arithmetical complexity of only ${\cal O}(m \log m \log (N/m)) $
and requires ${\cal O}(m \log (N/m))$ Fourier samples for this computation.
In contrast to other approaches there is no a priori knowledge needed about sparsity or  support bounds for  the vector ${\bf x}$. The algorithm automatically recognizes and exploits a possible short support of the vector and falls back to a usual radix-2 FFT algorithm if ${\bf x}$ has (almost) full support.
The numerical stability of the proposed algorithm is shown by numerical examples.}
\smallskip


\section{Introduction}

Algorithms for fast Fourier transform play a fundamental role in many areas in numerical analysis, particularly in signal and image processing. It is well-known that FFT algorithms for general vectors ${\bf x} \in {\C}^{N}$ require ${\cal O}(N \log N)$ arithmetical operations and that this qualitative bound cannot be improved, see \cite{morgen}. However, if the discrete Fourier transform is applied to recover vectors with special properties,
there is the hope for even faster algorithms.

In recent years, there has been some effort to derive new so-called ``sparse FFT'' algorithms that exploit the a priori knowledge that the vector to be recovered is sparse or has only a small amount of significant frequencies.
Often, further assumptions on the vector appear, as e.g.\ that the components to be recovered are from a certain  quantized range consisting of a finite set of real entries, see e.g.\ \cite{HIKP12a}.
Most of the proposed algorithms are based on randomization \cite{GIIS14, HIKP12a, HIKP12b, PR13} and achieve e.g.\ a complexity of ${\cal O}(k \log N)$ \cite{HIKP12a} for $k$-sparse signals, or even ${\cal O}(k \log k)$, see e.g.\ \cite{PR13}.
An obvious drawback of randomized approaches is that the algorithms do not always achieve the correct result (or an approximation of it) but only with a certain probability. Another problem is that there exists no sublinear algorithm to check the correctness of the result.

Completely deterministic sparse FFT algorithms for $k$-sparse signals have been proposed e.g.\ in \cite{Aka10,Aka14, I10, I13, LWC13}.
The underlying ideas are based on combinatorial approaches employing FFTs of different prime length and the Chinese remainder theorem. These algorithms usually have polynomial costs in $\log N$ and $k$ and only pay off for very large $N$ and strong sparsity.

Sparse FFT algorithms based on Prony's method, see \cite{HKPV13, PT14, PT15}, for $k$-sparse signals are usually based on singular value decompositions of size $k$ with a complexity of ${\cal O}(k^{3})$ and are therefore only efficient for small $k$. This complexity can be reduced in special cases using a splitting approach \cite{PT15}.

In a recent paper \cite{PW15}, the authors proposed a deterministic sparse FFT algorithm for vectors with short support that is based on usual FFT and is numerically stable with a complexity of ${\cal O}(m \log N)$ operations, where here $m$ denotes the support length of the signal.

\medskip

All algorithms mentioned above require a priori information as e.g.\ the exact sparsity or the support length of the vector to be recovered, or at least a suitable upper bound for it. They are just not applicable without this information or do not always achieve the correct recovery result by using only a guess for this bound.

However, in practice, while certain structures of the vector (as short support or sparsity) often appear, we do not always have the a priori knowledge on a good upper bound. Therefore it is of high interest to develop deterministic FFT algorithms that are able to automatically recognize certain structures of the vector during the algorithm and to exploit it suitably to reduce complexity and run time.

In this paper, we propose for the first time an algorithm that meets this requirement in the way that no a priori knowledge about the support length of the vector to be recovered is needed beforehand.
We present a new deterministic algorithm to recover a real non-negative vector $\mathbf x \in \R^N_{+}$ from its discrete Fourier transform $\widehat{\mathbf x}$. If  $\mathbf x$ has a support with support length $m$ being significantly  smaller than $N$, then the algorithm automatically recognizes this structure and provides the resulting vector with an arithmetical complexity of  $\mathcal{O}(m \log m \, \log\frac{N}{m})$ requiring at most ${\cal O}(m \log (N/m))$ Fourier data. 
The idea of the algorithm is based on divide-and-conquer techniques.

Direct applications for the reconstruction of sparse vectors from Fourier data, both with known or unknown support, appear for instance solving phase retrieval problems, 
where data have to be reconstructed from Fourier intensities. In this case, short support and positivity of the resulting vectors or images are frequently used preconditions in iterated projection algorithms, see e.g.\ \cite{fienup}.

This paper is structured as follows: After fixing the notations, we introduce the new fast algorithm in Section 2, together with detailed explanations of its structure and complexity. In Section 3, we apply the algorithm to noisy Fourier data and present some numerical results showing the numerical stability of the proposed algorithm in practice.

\subsection{Notations}
Let $\mathbf x \in \R^N_{+}$ with $N=2^{J}$ for some $J>0$ be a real vector with non-negative entries. We denote the discrete Fourier transform $\widehat{\mathbf x}$ of $\mathbf x$ by
$$ \widehat{\mathbf x} = \mathbf F_N \mathbf x, $$
where ${\mathbf F}_N := \left( \omega_N^{jk} \right)_{j,k=0}^{N-1} \in \C^{N\times N}$ is the Fourier matrix and $\omega_N := \e^{\frac{-2\pi{\mathrm i}}{N}}$.

The \emph{support length} $m = |\supp \mathbf x|$ of $\mathbf x \in \R^N_+$ is defined as the \emph{minimal positive integer} such that $x_k=0$ for $k \notin I:= \{(\mu+\ell)\mod N\ |\ \ell=0,\dots,m-1\}$. We call this index set $I$ the \emph{support index interval of $\mathbf x$}. The \emph{first support index of $\mathbf x$} is denoted by $\mu$.
Note that the first support index of ${\mathbf x}$ needs not to be the index of the first nonzero entry in ${\mathbf x}$.
Considering for example the vector ${\mathbf x} :=(13,21,0,0,0,10,31,0) \in {\R}_{+}^{8}$, we obtain a support length $m=5$, the support index interval $I=\{5,6,7,0,1 \}$ with the corresponding  signal values $(10,31,0,13,21)$,  and $\mu=5$, i.e., the support starts with $x_{5}=10$. Note further that the support index interval $I$ may contain indices corresponding to zero components of ${\bf x}$, as e.g.\ the index $7$ in the small example above. Therefore the support length $m$ is an upper bound of the sparsity, the number of nonzero entries of ${\mathbf x}$. But in any case, it holds that $x_\mu > 0$ and $x_{(\mu+m-1)\smod N} > 0$, i.e., the first and the last entry of the support of ${\bf x}$ are positive.

The support length of a vector $\mathbf x$ is always uniquely defined. However, the support index interval and the first support index $\mu$ are not necessarily unique. Consider e.g.\ the vector $\mathbf x \in \R^N_+$ with $x_0=x_{N/2}=1$ and $x_\ell=0$ for $\ell \in \{0,\dots,N-1\} \backslash \{0,N/2\}$. For this vector it is possible to choose either $\mu=0$ or $\mu=N/2$ whereas the support length is $N/2+1$ in both cases.

The \emph{periodized vectors} $\mathbf x^{(j)} \in \R^{2^j}_{+}$ of $\mathbf x$ are defined by
\begin{align}\label{periodization}
\mathbf x^{(j)} = (x_k^{(j)})_{k=0}^{2^j-1} := \left( \sum_{\ell=0}^{2^{J-j}-1} x_{k+2^j\ell} \right)_{k=0}^{2^j-1}
\end{align}
for $j=0,\dots,J$. In particular, $\mathbf x^{(0)} = \sum_{k=0}^{N-1} x_k$ is the sum of all components of $\mathbf x$, $\mathbf x^{(1)} = \left(\sum_{k=0}^{N/2-1} x_{2k}, \sum_{k=0}^{N/2-1} x_{2k+1}\right)^T$ and $\mathbf x^{(J)} = \mathbf x$.

We recall from  \cite{PW15},  that the components of the discrete Fourier transforms $\widehat{\mathbf x}^{(j)}$ need not to be computed but are already given as a subset of the set of components of $\widehat{\mathbf x} = \mathbf F_N \mathbf x$.

\begin{lemma}\label{lemma:periodization}
For the vectors $\mathbf x^{(j)} \in \R^{2^j}_+$, $j=0,\dots,J$, in (\ref{periodization}),  the discrete Fourier transform is given by
\[\widehat{\mathbf x}^{(j)} := \mathbf F_{2^j} \mathbf x^{(j)} = (\widehat{x}_{2^{J-j}k})_{k=0}^{2^j-1},\]
where $\widehat{\mathbf x} = (\widehat{\mathbf x}_k)_{k=0}^{N-1} = \mathbf F_N \mathbf x$ is the Fourier transform of $\mathbf x \in \R^N_+$.
\end{lemma}

\section{Sparse FFT algorithm for non-negative vectors}

Let us assume that the Fourier transform $\widehat{\mathbf x} = \mathbf F_N \mathbf x$ of $\mathbf x \in \R_{+}^N$ is given where $\mathbf x$ has only non-negative entries. We want to derive an algorithm that automatically recognizes  a possible shorter support of ${\bf x}$ and applies in this case a faster
reconstruction algorithm while falling back to the usual FFT with ${\cal O}(N \log_{2} N)$ complexity if ${\bf x}$ has full support or does not possess a support with a length being significantly smaller than $N$.

The main idea  to reconstruct ${\bf x}$ from $\widehat{\bf x}$ is to employ the divide-and-conquer technique  similarly as for usual radix-2 FFT.
Starting with   ${\bf  x}^{(0)} = \widehat{x}_{0} \in {\R}_{+}$  we recover ${\bf x}$ iteratively  by
 reconstructing $\mathbf x^{(j+1)}$ from $\mathbf x^{(j)}$ for $j=0,\dots,J-1$. At each level, we check the support length $m_{j} := |\supp \mathbf x^{(j)}|$ of the vector $\mathbf x^{(j)}$ (that has been computed at the previous level) and distinguish two cases: $m_{j} > 2^{j-1}$ and $m_{j} \le 2^{j-1}$. 
In the first case, i.e., if the support length of $\mathbf x^{(j)}$ is greater than half of the vector length, we cannot benefit from a short support and therefore compute $\mathbf x^{(j+1)}$ using ${\bf x}^{(j)}$ and employing an FFT algorithm of length $2^{j}$.
 In the second case, if the support of $\mathbf x^{(j)}$ is shorter than half of the vector length, we apply a modified  reconstruction algorithm that only requires ${\cal O}(m_{j} \log m_{j})$ floating point operations.

In both cases, we aim at computing $\mathbf x^{(j+1)}$ by a simplified inverse fast Fourier transform that exploits the values of $\mathbf x^{(j)}$ being known already from the previous iteration step. 

Splitting the vector ${\bf x}^{(j+1)}$ into the two partial vectors $\mathbf x_0^{(j+1)} = (x_k^{(j+1)})_{k=0}^{2^j-1}$ and $\mathbf x_1^{(j+1)} = (x_k^{(j+1)})_{k=2^j}^{2^{j+1}-1}$ of length $2^{j}$, we recall that by definition 
\begin{equation}\label{e1}
\mathbf x^{(j)} = \mathbf x_0^{(j+1)} + \mathbf x_1^{(j+1)}.
\end{equation} 
Using this equation and Lemma \ref{lemma:periodization} we observe the following relation for  $\widehat{\mathbf x}^{(j+1)}$,
\begin{align*}
\widehat{\mathbf x}^{(j+1)} &= (\widehat{x}_{2^{J-j-1}k})_{k=0}^{2^{j+1}-1}  = \mathbf F_{2^{j+1}} \mathbf x^{(j+1)} = (\omega_{2^{j+1}}^{k\ell})_{k,\ell=0}^{2^{j+1}-1} \begin{pmatrix} \mathbf x_0^{(j+1)}\\\mathbf x_1^{(j+1)} \end{pmatrix}\\
	&= (\omega_{2^{j+1}}^{k\ell})_{k,\ell=0}^{2^{j+1}-1,2^j-1} \mathbf x_0^{(j+1)} + (\omega_{2^{j+1}}^{k\ell})_{k=0,\ell=2^j}^{2^{j+1}-1,2^{j+1}-1} (\mathbf x^{(j)} - \mathbf x_0^{(j+1)})\\
	&= (\omega_{2^{j+1}}^{k\ell})_{k,\ell=0}^{2^{j+1}-1,2^j-1} \mathbf x_0^{(j+1)} + ((-1)^k \omega_{2^{j+1}}^{k\ell})_{k,\ell=0}^{2^{j+1}-1,2^j-1} (\mathbf x^{(j)} - \mathbf x_0^{(j+1)}). \\
	\end{align*}
While the even-indexed Fourier components of $\widehat{\mathbf x}^{(j+1)}$ only contain information on $\mathbf x^{(j)}$
by Lemma \ref{lemma:periodization}, the odd components give new information on the vector ${\bf x}^{(j+1)}$. Restricting the equation system to the odd components yields
\begin{align}\label{equation:fouriervalues}
(\widehat{x}^{(j+1)}_{2k+1})_{k=0}^{2^{j}-1} &= (\widehat{x}_{2^{J-j-1}(2k+1)})_{k=0}^{2^j-1} \nonumber \\
& = (\omega_{2^{j+1}}^{(2k+1)\ell})_{k,\ell=0}^{2^j-1} \mathbf x_0^{(j+1)} - (\omega_{2^{j+1}}^{(2k+1)\ell})_{k,\ell=0}^{2^j-1} (\mathbf x^{(j)} - \mathbf x_0^{(j+1)})\nonumber\\
	&=  (\omega_{2^{j+1}}^{(2k+1)\ell})_{k,\ell=0}^{2^j-1} (2 \mathbf x_0^{(j+1)} -  \mathbf x^{(j)})\\
	&=  \mathbf F_{2^j} \cdot\mathrm{diag}(\omega_{2^{j+1}}^\ell)_{\ell=0}^{2^j-1} \, (2 \mathbf x_0^{(j+1)} - \mathbf x^{(j)}).\nonumber
\end{align}

\bigskip
\noindent
\textbf{(1) First case:} $m_{j} > 2^{j-1}$\\
In this case we just exploit the observations in (\ref{equation:fouriervalues}) and obtain
\begin{align*}
\mathbf x_0^{(j+1)} = \frac{1}{2} \left( \mathrm{diag}(\omega_{2^{j+1}}^{-\ell})_{\ell=0}^{2^j-1} \cdot \mathbf F_{2^j}^{-1} \cdot (\widehat{x}_{2^{J-j-1}(2k+1)})_{k=0}^{2^j-1} + \mathbf x^{(j)} \right).
\end{align*}
Thus, ${\bf x}_{0}^{(j+1)}$ can be computed via an inverse FFT of length $2^{j}$. Further we require $2^{j}$ complex multiplications, $2^{j}$ additions and one dyadic shift by $2$. Finally,  ${\bf x}_{1}^{(j+1)} = {\bf x}^{(j)} - {\bf x}_{0}^{(j+1)}$ is obtained by ${\cal O}(2^{j})$ flops.

\bigskip
\noindent
\textbf{(2) Second case:} $m_{j} \le 2^{j-1}$\\
For $m_{j} \le 2^{j-1}$, we first compute $L_{j} := \lceil \log_2 m_{j} \rceil \le {j-1}$. 
We denote the first support index of $\mathbf x^{(j)}$ by $\mu^{(j)}$.
Then the support of $\mathbf x^{(j)}$ lies within an interval of length $2^{L_{j}}$, beginning at the index $\mu^{(j)}$. 
Instead of ${\bf x}^{(j)}, \, {\bf x}^{(j+1)}_{0}$, and ${\bf x}^{(j+1)}_{1}$ of length $2^{j}$ as above, we consider now the partial vectors of length $2^{L_{j}}$ containing the relevant support  
\begin{align*}
\widetilde{\mathbf x}^{(j)} & = (x_{(\mu^{(j)}+r)\smod 2^j}^{(j)})_{r=0}^{2^{L_{j}}-1}, \\
\widetilde{\mathbf x}_0^{(j+1)} & = (x_{(\mu^{(j)}+r)\smod 2^j}^{(j+1)})_{r=0}^{2^{L_{j}}-1}, \quad \widetilde{\mathbf x}_1^{(j+1)} = (x_{2^j+(\mu^{(j)}+r)\smod 2^j}^{(j+1)})_{r=0}^{2^{L_{j}}-1}.
\end{align*} 

Indeed, by (\ref{e1}) the vector $\mathbf x^{(j+1)}$ that we want to reconstruct cannot have more than $2^{L_{j}+1}$ positive entries, and these relevant entries are contained in the restricted vectors $\widetilde{\mathbf x}_0^{(j+1)}$ and $\widetilde{\mathbf x}_1^{(j+1)}$. Since the condition 
$$ \widetilde{\mathbf x}_0^{(j+1)}+ \widetilde{\mathbf x}_1^{(j+1)} = \widetilde{\mathbf x}^{(j)} $$
is still satisfied, we only need  $2^{L_{j}}$ further linearly independent conditions to recover $\mathbf x^{(j+1)}$ completely.
Employing 
the equation system (\ref{equation:fouriervalues}) and using the shorter support of the partial vectors, we find
\begin{align}\label{equation:fouriervaluesII}
(\widehat{x}_{2^{J-j-1}(2k+1)})_{k=0}^{2^j-1} &=  (\omega_{2^{j+1}}^{(2k+1)\ell})_{k,\ell=0}^{2^j-1} (2 \mathbf x_0^{(j+1)} - \mathbf x^{(j)})\nonumber\\
	&= (\omega_{2^{j+1}}^{(2k+1)((\mu^{(j)}+r) \smod 2^j)})_{k,r=0}^{2^j-1, 2^{L_{j}}-1} (2 \widetilde{\mathbf x}_0^{(j+1)}
- \widetilde{\mathbf x}^{(j)}).\nonumber
\end{align}

Instead of considering these $2^{j}$ equations for $k=0, \ldots , 2^{j}-1$, we employ now only the $2^{L_{j}}$ equations for $k= 2^{j-L_{j}}p$, $p=0, \ldots , 2^{L_{j}}-1$, and obtain 
\begin{eqnarray} \label{case2}
(\widehat{x}_{2^{J-L_{j}}p+2^{J-j-1}})_{p=0}^{2^{L_{j}}-1} &= &
(\widehat{x}_{2^{J-j-1}(2^{j+1-L_{j}}p+1)})_{p=0}^{2^{L_{j}}-1} \\
& = & (\omega_{2^{j+1}}^{(2^{j+1-L_{j}}p+1)((\mu^{(j)}+r) \smod 2^j)})_{p,r=0}^{2^{L_{j}}-1, 2^{L_{j}}-1}  (2 \widetilde{\mathbf x}_0^{(j+1)}
- \widetilde{\mathbf x}^{(j)}), \nonumber
\end{eqnarray}
where 
\begin{align*}
& (\omega_{2^{j+1}}^{(2^{j+1-L_{j}}p+1)((\mu^{(j)}+r) \smod 2^j)})_{p,r=0}^{2^{L_{j}}-1, 2^{L_{j}}-1} \\
&= (\omega_{2^{L_{j}}}^{p(\mu^{(j)}+r)})_{p,r=0}^{2^{L_{j}}-1} \, \mathrm{diag} (\omega_{2^{j+1}}^{(\mu^{(j)}+r) \smod 2^{j}})_{r=0}^{2^{L_{j}}-1}\\
&= \mathrm{diag} (\omega_{2^{L_{j}}}^{\mu^{(j)}p})_{p=0}^{2^{L_{j}}-1} \, {\bf F}_{2^{L_{j}}} \, \mathrm{diag} (\omega_{2^{j+1}}^{(\mu^{(j)}+r) \smod 2^{j}})_{r=0}^{2^{L_{j}}-1}.
\end{align*}
We finally conclude from (\ref{case2})
\begin{eqnarray*} \widetilde{\mathbf x}_0^{(j+1)} &=& 
\frac{1}{2}
 \mathrm{diag} (\omega_{2^{j+1}}^{-(\mu^{(j)}+r) \smod 2^{j}})_{r=0}^{2^{L_{j}}-1} \, {\bf F}_{2^{L_{j}}}^{-1} \, \mathrm{diag} (\omega_{2^{L_{j}}}^{-\mu^{(j)}p})_{p=0}^{2^{L_{j}}-1} \,
(\widehat{x}_{2^{J-L_{j}}p+2^{J-j-1})})_{p=0}^{L_{j}-1}\\
& &  + \frac{1}{2} \widetilde{\mathbf x}^{(j)}.
\end{eqnarray*}
Thus, in this case the recovery of $\widetilde{\mathbf x}_0^{(j+1)}$ requires only ${\cal O}(2^{L_{j}}  L_{j}) = {\cal O}(m_{j} \log_{2} m_{j})$ flops, while $\widetilde{\mathbf x}_1^{(j+1)}$ is obtained from $ \widetilde{\mathbf x}_1^{(j+1)} = \widetilde{\mathbf x}^{(j)} - \widetilde{\mathbf x}_0^{(j+1)}$.

The new fast algorithm to compute the vector $\mathbf x \in \R^N_{+}$ with possible short support from its Fourier transform $\widehat{\bf x}$ can be summarized as in the following algorithm, where we iteratively compute the periodized vectors $\mathbf x^{(j)} \in {\R}^{2^{j}}_+$ by applying either the method of case 1 or case 2 at each iteration level.
If there is some a priori information available on a lower bound $2^{s-1}$ for the support length $m$ of ${\mathbf x}$, then we may start the iteration by computing the periodized vector ${\bf x}^{(s)}$ of length $2^{s}$, otherwise we just start with $s=0$ in the algorithm.

\begin{algorithm}\null \label{algorithm} (Sparse FFT for real non-negative vectors)\\
\textbf{Input:} $\widehat{\mathbf x} = (\widehat{x}_{k})_{k=0}^{N-1}\in \C^N$, $N=2^J$; \\
\phantom{\textbf{Input:}} $s=0$ or $s$ such that $2^{s-1}$ is a lower bound for $ m =|{\rm supp} \, {\mathbf x}|$; \\
\phantom{\textbf{Input:}} threshold parameter $T$.
\begin{enumerate}
\item Generate $\widehat{\mathbf x}^{(s)} := (\widehat{x}_{2^{J-s}k})_{k=0}^{2^s-1}$ by extracting suitable components from $\widehat{\mathbf x}$.
\item Compute the periodized vector $\mathbf x^{(s)} := \mathbf F_{2^s}^{-1} \widehat{\mathbf x}^{(s)}$ by inverse FFT of length $2^s$.
\item For $k=0, \ldots , 2^{s}-1$, apply a threshold procedure
$$ x_{k}^{(s)} := \left\{ \begin{array}{ll} {\rm Re} \, x_{k}^{(s)} & {\it  if} \, {\rm Re} \, x_{k}^{(s)} \ge T, \\
	0 & {\it else.} \end{array} \right. $$
\item For $j=s,\dots,J-1$ do \\
	Compute $m_{j} := |\supp \mathbf x^{(j)}|$ and find the first support index $\mu^{(j)}$ of $\mathbf x^{(j)}$.
	\begin{itemize}
	\item {\bf\rm Case 1}: If $m_{j} > 2^{j-1}$, then \\
	Build $ {\mathbf y}^{(j)} := (\widehat{x}_{2^{J-j-1}(2k+1)})_{k=0}^{2^j-1}$
	and compute
	\[\mathbf z^{(j)} := \mathrm{diag}(\omega_{2^{j+1}}^{-\ell})_{\ell=0}^{2^j-1} \cdot \mathbf F_{2^j}^{-1} \cdot \mathbf y^{(j)} \]
	using an inverse FFT of length $2^{j}$.\\
	Compute 
	\[{\mathbf x}^{(j+1)} := \frac{1}{2} \begin{pmatrix} \mathbf x^{(j)} + \mathbf z^{(j)}\\ \mathbf x^{(j)} - \mathbf z^{(j)} \end{pmatrix}.\]
	For $k=0, \ldots , 2^{j+1}-1$, apply a threshold procedure
	$$ x_{k}^{(j+1)} := \left\{ \begin{array}{ll} {\rm Re} \, x_{k}^{(j+1)} & {\it  if} \, {\rm Re} \, x_{k}^{(j+1)} \ge T, \\
	0 & {\it else.} \end{array} \right. $$
	end (if).
	\item {\bf\rm Case 2}: If $m_{j} \le 2^{j-1}$, then \\
	Compute $L_{j} := \lceil \log_2 m_{j} \rceil$ and $\mu^{(j)}$.
	Build the vectors
	\[\widetilde{\mathbf x}^{(j)} := (x_{(\mu^{(j)}+r)\smod 2^j}^{(j)})_{r=0}^{2^{L_{j}}-1}, \quad \mathbf y^{(j)} := (\widehat{x}_{2^{J-L_{j}}p+2^{J-j-1}})_{p=0}^{2^{L_{j}}-1}.\]
	Compute
	$$ {\mathbf z}^{(j)} := \mathrm{diag} (\omega_{2^{L_{j}}}^{-\mu^{(j)}p})_{p=0}^{2^{L_{j}}-1} \, {\bf F}_{2^{L_{j}}}^{-1} \, \mathrm{diag} (\omega_{2^{j+1}}^{-(\mu^{(j)}+r) \smod 2^{j}})_{r=0}^{2^{L_{j}}-1} \,
{\mathbf y}^{(j)}$$
using an inverse FFT of length $2^{L_{j}}$.\\

	Compute $\widetilde{\mathbf x}_0^{(j+1)} := \frac{1}{2}(\widetilde{\mathbf x}^{(j)} + \mathbf z^{(j)})$ and $\widetilde{\mathbf x}_1^{(j+1)} := \frac{1}{2}(\widetilde{\mathbf x}^{(j)} - \mathbf z^{(j)})$.\\
	For $k=0, \ldots , 2^{L_{j}}-1$, apply a threshold procedure
	\begin{eqnarray*} 
	(\widetilde{\mathbf x}_{0}^{(j+1)})_{k} &:=& \left\{ \begin{array}{ll} {\rm Re} \, (\widetilde{\mathbf x}_{0}^{(j+1)})_{k} & {\it  if} \, {\rm Re} \, (\widetilde{\mathbf x}_{0}^{(j+1)})_{k} \ge T, \\
	0 & {\it else.} \end{array} \right. \\
	(\widetilde{\mathbf x}_{1}^{(j+1)})_{k} &:=& \left\{ \begin{array}{ll} {\rm Re} \, (\widetilde{\mathbf x}_{1}^{(j+1)})_{k} & {\it  if} \, {\rm Re} \, (\widetilde{\mathbf x}_{1}^{(j+1)})_{k} \ge T, \\
	0 & {\it else.} \end{array} \right.
	\end{eqnarray*}
	Determine $\mathbf x_0^{(j+1)}$ and $\mathbf x_1^{(j+1)}$ by
	\[(\mathbf x_0^{(j+1)})_{(\mu^{(j)}+k) \smod 2^j} := \left\{ \begin{array}{ll} (\widetilde{\mathbf x}_0^{(j+1)})_k & \qquad k=0, \dots, 2^{L_{j}}-1,\\ 0 & \qquad k=2^{L_{j}}, \dots, 2^j, \end{array} \right. \]
	\[(\mathbf x_1^{(j+1)})_{(\mu^{(j)}+k) \smod 2^j} := \left\{ \begin{array}{ll} (\widetilde{\mathbf x}_1^{(j+1)})_k & \qquad k=0, \dots, 2^{L_{j}}-1,\\ 0 & \qquad k=2^{L_{j}}, \dots, 2^j. \end{array} \right.\]
	Set $\mathbf x^{(j+1)} := \begin{pmatrix} \mathbf x_0^{(j+1)} \\ \mathbf x_1^{(j+1)} \end{pmatrix}$. \\
	end (if)
\end{itemize}

end (for)
\end{enumerate}
\textbf{Output:} $\mathbf x^{(J)} = \mathbf x$.
\end{algorithm}

For a Matlab implementation of this algorithm we refer to our homepage \newline
{\small{\tt  http://na.math.uni-goettingen.de}}.

The threshold parameter $T$ in the algorithm ensures that we obtain a real non-negative vector as a result regardless of small numerical errors that may arise. In case of noisy Fourier data, the threshold parameter needs to be chosen suitably to suppress errors in the solution vector, see Section \ref{res}.

\begin{remark} 
1. At every reconstruction step the algorithm automatically decides whether the first or the second case applies. For this purpose, the support length of $\mathbf x^{(j)}$ has to be computed. This can be efficiently done by using the known support indices of the preceding periodization $\mathbf x^{(j-1)}$. By definition of the periodization, the vector $\mathbf x^{(j)}$ can only have positive entries at the support indices of $\mathbf x^{(j-1)}$ and at these indices shifted by $2^{j-1}$. Hence, only $2m_{j-1}$ entries have to be considered in order to find the support length $m_{j}$ and the first support index $\mu^{(j)}$ of $\mathbf{x}^{(j)}$ causing an effort of ${\cal O}(m_{j})$ flops. 

\smallskip

2. We want to emphasize the importance of the modulo operation in the diagonal matrix $\mathrm{diag} (\omega_{2^{j+1}}^{-(\mu^{(j)}+r) \smod 2^{j}})_{r=0}^{2^{L_{j}}-1}$ in case 2 of the algorithm. While for $\mu^{(j)} \le 2^{j} - 2^{L_{j}}$ this matrix can be simplified to 
$$ \mathrm{diag} (\omega_{2^{j+1}}^{-(\mu^{(j)}+r) \smod 2^{j}})_{r=0}^{2^{L_{j}}-1} = \mathrm{diag} (\omega_{2^{j+1}}^{-(\mu^{(j)}+r) })_{r=0}^{2^{L_{j}}-1} =  \omega_{2^{j+1}}^{-\mu^{(j)}} \, \mathrm{diag} (\omega_{2^{j+1}}^{-r})_{r=0}^{2^{L_{j}}-1}, $$
 it follows for $\mu^{(j)} > 2^{j} - 2^{L_{j}}$ that 
$$ \mathrm{diag} (\omega_{2^{j+1}}^{-(\mu^{(j)}+r) \smod 2^{j}})_{r=0}^{2^{L_{j}}-1} = 
\omega_{2^{j+1}}^{-\mu^{(j)}} \, \left(
\begin{array}{cc}  {\bf I}_{2^{j}-\mu^{(j)}} & \\  & -{\bf I}_{2^{L_{j}}- 2^{j} + \mu^{(j)}} \end{array} \right) \, \mathrm{diag} (\omega_{2^{j+1}}^{-r})_{r=0}^{2^{L_{j}}-1} ,
$$
where ${\bf I}_{2^{j}-\mu^{(j)}}$ and ${\bf I}_{2^{L_{j}}- 2^{j} + \mu^{(j)}}$ denote identity matrices of the given size.
\end{remark}

Let us summarize the numerical effort of the complete algorithm.
If we start with $s=0$, the complexity of the algorithm is at most $\mathcal{O}(m \log_2 m \log_2\frac{N}{m})$. 
At each iteration step, the support length $m_{j}$ of the periodized vectors $\mathbf x^{(j)}$ can only increase, 
i.e., we have $m_{0} \le m_{1} \le \ldots \le m_{J-1} \le m_{J}=m$.
Once the final support length has been achieved, say at the iteration step $L$, and $2^{L-1} < m_{L} = m \le 2^{L}$, we will always employ the 
second case in the further iteration steps $j=L+1, \ldots, J-1$ that 
requires ${\cal O}(m \log_{2} m)$ flops at each step.

The first $L$ reconstructing steps $j=1, \ldots , L$ may require either the first or the second case (depending on the distribution of nonzero values of ${\bf x}$) and require at most $\mathcal{O}(2^{L} L ) = \mathcal{O}(m \log m )$ flops, caused by the inverse FFT of size $2^j$, a multiplication with a diagonal matrix of size $2^j \times 2^j$ and $2^{j+2}$ additions and multiplications computing the periodization $\mathbf x^{(j+1)}$ at each iteration step, similarly as a usual FFT algorithm of length $2^{L}$. 

Together, we thus require ${\cal O}((J-L) m \log_{2} m) = {\cal O}(m \log_{2} m (\log_{2} N/m))$ flops to compute $\mathbf x$.

\begin{remark}
1. The sublinear complexity of the algorithm can only be achieved by employing less than the given $N$ Fourier samples in the vector $\widehat{\mathbf x}$. Indeed at the $j$-th iteration step, we use either $2^{j}$ new Fourier samples in the first case or only $2^{L_{j}}$ Fourier samples in the second case, collected in the vector   ${\mathbf y}^{(j)}$. Assuming as before that $2^{L-1} < m \le 2^{L}$, we apply at steps $L+1, \ldots , J-1$ the second case requiring only $(J-L-1) m$ Fourier samples while we need  at most $2^{L+1}$ Fourier samples at the first
steps $1, \ldots , L$. Altogether, the number of applied Fourier samples is bounded by ${\cal O}( m (\log_{2}N/m ))$. 

2. The proposed algorithm is efficient  for any vector $\mathbf x \in \R^N_{+}$, whether or not it has short support. The complexity does never exceed  $\mathcal{O}(N \log_{2} N)$ of usual FFT algorithms. In case that the support of the vector is quite short compared to the full vector length, we benefit from the algorithm concerning the computational complexity.
We may even benefit from the algorithm if the vector has (almost) full support  length but is sparse, such that 
the second case applies in intermediate steps. 
For example, a vector ${\bf x} \in {\R}^{N}_{+}$ containing  several equidistantly distributed short support pieces, these support pieces may add up to one short support interval for smaller periodized vectors ${\bf x}^{(j)}$  such that we can take advantage of the algorithm.

Let us give an example: Choose $\mathbf x \in \R^{1024}_+$ with positive entries $x_0=1$, $x_{256}=1$, $x_{512}=1$ and $x_{768}=1$. Then $\mathbf x^{(9)}$ has two positive entries: $x_0=2$ and $x_{256}=2$. All further periodizations $\mathbf x^{(8)}, \dots, \mathbf x^{(0)}$ only have one positive entry: $x_0=4$ such that 
case 2 applies with $L_{j} =0$ for j=1, \ldots , 8.
\end{remark}

\section{Numerical Results}
\label{res}

We consider the numerical stability of the proposed algorithm. 
For that purpose, we apply the algorithm to Fourier data being perturbed  by uniform noise $\boldsymbol{\varepsilon} = (\varepsilon_k)_{k=0}^{N-1}$, i.e., we have given data
\[\widehat{y}_k = \widehat{x}_k + \varepsilon_k\]
with $|\varepsilon_k| \le \delta$. The above algorithm also applies to noisy data where we have to modify the threshold parameter $T$ suitably. 

In the noisy case, it is of particular importance to determine the support index interval correctly at each iteration step. As before, we do this by only considering the relevant entries given by the support index interval of the preceding periodization. Additionally, the threshold parameter $T$ has to be set in order to distinguish between relevant components of ${\bf x}$ and noise.

\begin{figure}[t!]
\begin{center}
\includegraphics[width=0.9\linewidth]{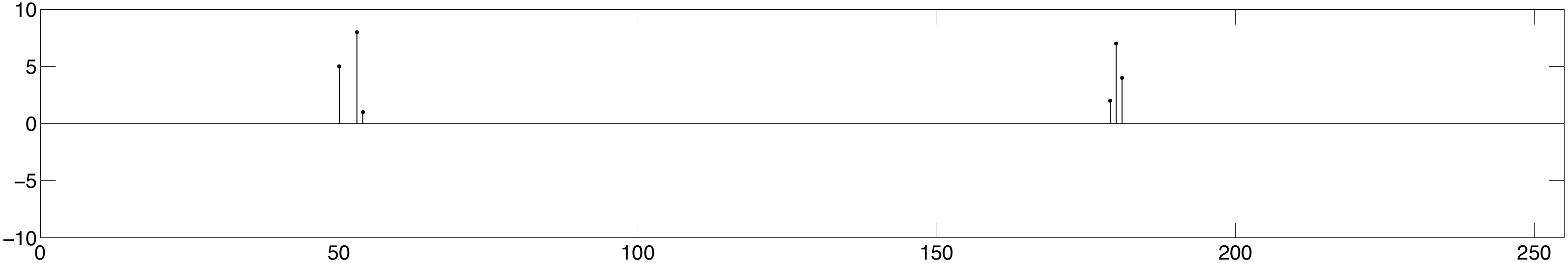}\\
{\small (a)} \\[1ex]
\includegraphics[width=0.9\linewidth]{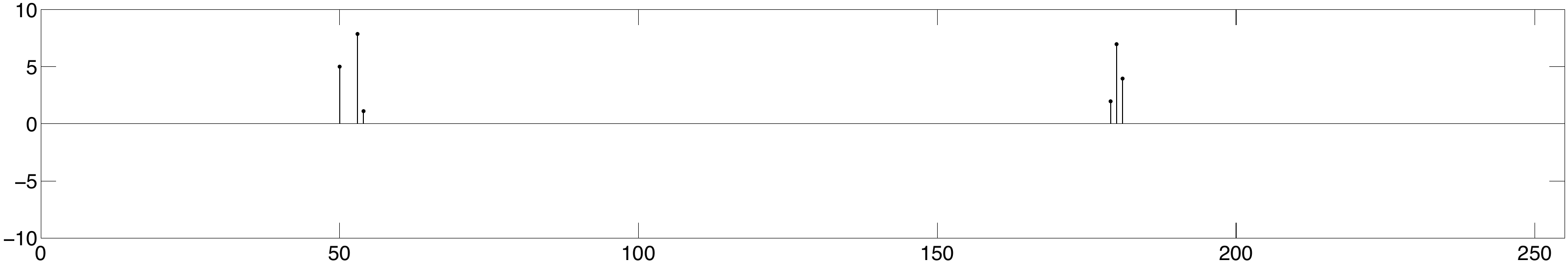}\\
{\small (b)}\\[1ex]
\includegraphics[width=0.9\linewidth]{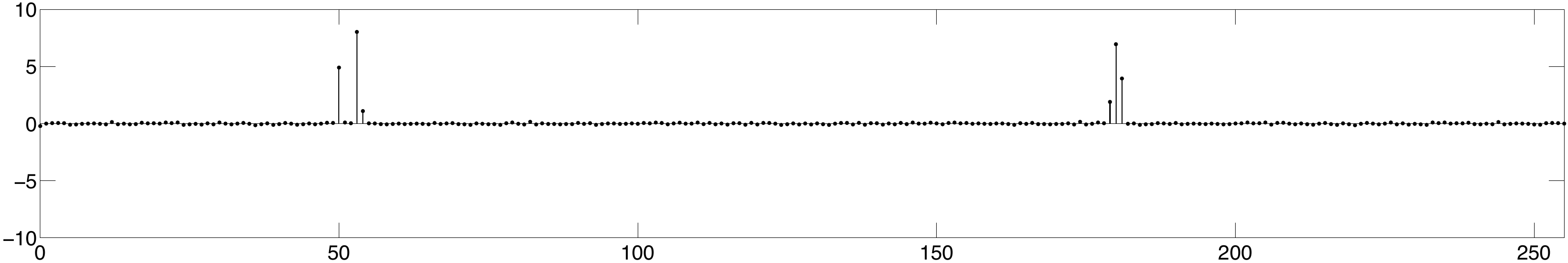}\\
{\small (c)}\\[1ex]
\caption{(a) Original vector $\mathbf x$ of length $N=256$; (b) Reconstruction of $\mathbf x$ using the sparse FFT Algorithm \ref{algorithm}; (c) Reconstruction of $\mathbf x$ using the inverse FFT.}
\label{example}
\end{center}
\end{figure}

Let us now give some numerical examples. We measure the noise level for the Fourier data using the signal-to-noise-ratio
\[\text{SNR} = 20 \cdot \log_{10}\frac{\|\widehat{\mathbf x}\|_2}{\|\boldsymbol{\varepsilon}\|_2}.\]
The error of the reconstruction is given by $\|\mathbf x - \mathbf x'\|_2/N$, where $\mathbf x'$ denotes the reconstruction of $\mathbf x$ by our algorithm.

Consider first a vector $\mathbf x$ of length $N=2^8=256$, with nonzero entries $x_{50}=5$, $x_{53}=8$, $x_{54}=1$, $x_{179}=2$, $x_{180}=7$ and $x_{181}=4$. We disturb the Fourier data $\widehat{\mathbf x}$ by uniform noise $\boldsymbol{\varepsilon}$ with $\text{SNR} = 20$ and reconstruct $\mathbf x$ from $\widehat{\mathbf y} = \widehat{\mathbf x} + \boldsymbol{\varepsilon}$ using our algorithm. In this example, we have $\|\boldsymbol{\varepsilon}\|_\infty = 2.149$ and $\|\boldsymbol{\varepsilon}\|_1/N = 1.100$. Choosing $T = 0.9$, the algorithm performs seven reconstruction steps in order to recover $\mathbf x$, where in the first three steps, case 1 is used and in the last four steps, case 2 applies.

The reconstructed vector $\mathbf x'$ has nonzero entries $x'_{50}=5.005$, $x'_{53}=7.868$, $x'_{54}=1.111$,  $x'_{179}=1.972$, $x'_{180}=6.980$, and $x'_{181}=3.964$ yielding an error $\|\mathbf x - \mathbf x'\|_2/N = 7.003 \cdot 10^{-4}$. In particular, the support of ${\bf x}$ is correctly found.
Compared to this, the inverse Fourier transform provides  an error $\|\mathbf x - \mathbf F_{256}^{-1} \widehat{\mathbf y}\|_2/N = 0.0049$.

Figure \ref{example} illustrates the vector $\mathbf x$ as well as both reconstructions from $\widehat{\mathbf y}$, by our algorithm and by an inverse FFT.

\begin{figure}[t!]
\begin{center}
\includegraphics[width=0.8\linewidth]{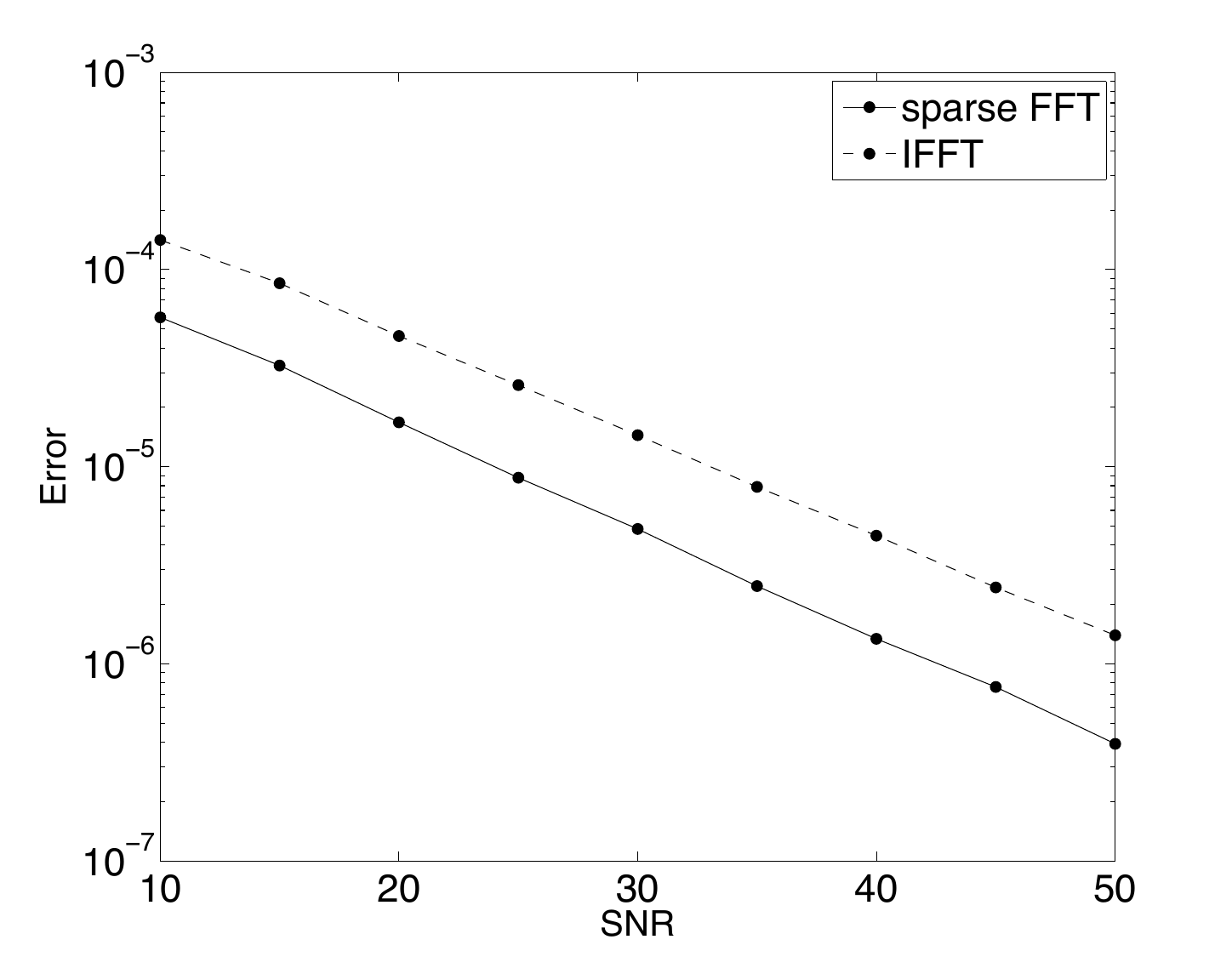}\\
\caption{Average reconstruction error $\|\mathbf x - \mathbf x'\|_2/N$ for different levels of uniform noise, comparing our deterministic sparse FFT algorithm and usual inverse FFT.}
\label{alg_vs_ifft}
\end{center}
\end{figure}

In a second example, we apply the algorithm to reconstruct randomly chosen real vectors $\mathbf x \in \R^N_{+}$ for $N=2^{15}$ with entries $0 \le x_k \le 10$ from their Fourier data and analyze the reconstruction error. The vectors in our experiment have a support length of $m = 15$. For each noise level between $\text{SNR} = 10$ and $\text{SNR} = 50$, we consider 100 randomly chosen vectors and apply the algorithm to the noisy Fourier data $\widehat{\mathbf y}$. The reconstruction error $\|\mathbf x - \mathbf x'\|_2/N$ for a reconstruction $\mathbf x'$ by our algorithm is computed as well as the error $\|\mathbf x - \mathbf F_{N}^{-1} \widehat{\mathbf y}\|_2/N$ of an inverse FFT reconstruction. The parameter $T$ is chosen appropriately for each noise level. The results of the experiment are shown in Figure \ref{alg_vs_ifft}.

The findings show that the proposed algorithm is numerically stable and returns a good reconstruction with a small error compared to usual inverse FFT. Additionally, it has a lower complexity: in our example with $N=2^{15}$, our Algorithm \ref{algorithm} applies in average $4.5$ times case 1 and  $10.5$ times case 2.

\subsection*{Acknowledgement}

The research in this paper is funded by the project PL 170/16-1 of the German Research Foundation (DFG). This is gratefully acknowledged.


\begin{thebibliography}{18}

\bibitem{Aka10}
A. Akavia, Deterministic sparse Fourier approximation  via fooling arithmetic progressions,
in Proc. 23rd COLT, 2010, pp. 381--393.

\bibitem{Aka14}
A. Akavia, Deterministic sparse Fourier approximation  via approximating arithmetic progressions,
IEEE Trans. Inform. Theory {\bf 60}(3) (2014), 1733--1741.

\bibitem{fienup}
J.~R. Fienup, {{Phase retrieval algorithms: A comparison},} Appl.
  Opt. \textbf{21}, 2758--2769 (1982).
  
\bibitem{GIIS14}
A. Gilbert, P. Indyk, M.A. Iwen, and L. Schmidt,
Recent developments in the sparse Fourier transform,
IEEE Signal Processing Magazine {\bf 31}(5) (2014), 91--100.


\bibitem{HIKP12a}
H. Hassanieh, P. Indyk, D. Katabi, and E. Price, Nearly optimal algorithm  for  sparse Fourier transform, Proc. 44th annual ACM symposium on Theory of Computing, 2012, 
pp. 563--578.


\bibitem{HIKP12b}
H. Hassanieh, P. Indyk, D. Katabi, and E. Price, Simple and practical algorithm for 
sparse Fourier transform,  Proc. 23th Annual ACM-SIAM Symposium on Discrete Algorithms (SODA '12), 2012, 
pp. 1183--1194.

\bibitem{HKPV13}
S. Heider, S. Kunis, D. Potts, and M. Veit, A sparse Prony FFT, Proc. 10th International Conference on Sampling Theory and Applications (SAMPTA), 2013, pp. 572--575.



\bibitem{I10} M.A. Iwen, Combinatorial sublinear-time Fourier algorithms,  Found. Comput. Math. {\bf 10} (2010), 303--338.

\bibitem{I13}
M.A.  Iwen, Improved approximation guarantees for sublinear-time Fourier algorithms, 
Appl. Comput. Harmon. Anal. {\bf 34} (2013), 57--82.

\bibitem{LWC13}
D. Lawlor, Y. Wang, and A. Christlieb,  Adaptive sub-linear time Fourier algorithms, Advances in Adaptive Data Analysis {\bf 5}(1)
(2013), 1350003 (25 pages).

\bibitem{morgen}
J. Morgenstern,  Note on a lower bound on the  linear complexity of the fast Fourier transform, J. Assoc. Comput. Mach. {\bf 20}(2) (1973), 305--306.

\bibitem{PR13}
S. Pawar and K. Ramchandran, Computing  a $k$-sparse $n$-length discrete Fourier transform  using at most  $4k$ samples  and ${\cal O}(k \log k)$ complexity, IEEE International Symposium on Information Theory (2013),  pp. 464--468.

\bibitem{PT14}
G. Plonka and M. Tasche, Prony methods for recovery of structured functions,
GAMM-Mitt. {\bf 37}(2) (2014),  239--258.

\bibitem{PW15}
G. Plonka and K. Wannenwetsch, A deterministic sparse FFT algorithm for vectors with small support, Numer. Algor. {\bf 71}(4) (2016), 889--905.

\bibitem{PT15}
D. Potts, M. Tasche, T. Volkmer, Efficient spectral estimation by MUSIC and ESPRIT with application to sparse FFT,  Front. Appl. Math. Stat., 29 February 2016.

%
%
\end{thebibliography}

\end{document}